\documentclass[12pt]{article}

\usepackage{graphicx} 
\usepackage{amsmath}

\newtheorem{theorem}{Theorem}

\newcommand\BV{Brunt-V{\"a}is{\"a}l{\"a}\ }

\newcommand\pde{pde\ }

\newcommand\sfun{\psi}

\title{On functional equations leading to exact solutions for standing internal waves}
\author{F. Beckebanze (1) 
\ and
G. Keady (2) 
\thanks{Email address for correspondence Grant.Keady@curtin.edu.au}\\
( (1) Institute for Marine and Atmospheric Research Utrecht, \\
 \textsc{The Netherlands} \\
(2) Department of Mathematics, Curtin University,\\
 \textsc{Australia} ) }

\date{\today} 

\begin{document}

\maketitle

\begin{abstract}
The Dirichlet problem for the wave equation is a classical example of a problem which is not well-posed.
Nevertheless, it has been used to model internal waves oscillating sinusoidally in time,  in various situations, standing internal waves amongst them.
We consider internal waves in two-dimensional domains bounded above by the plane $z=0$ and below by $z=-d(x)$ for depth functions $d$. 
This paper draws attention to the Abel and Schr{\"o}der functional equations as  a convenient way of organizing analytical solutions. 
Exact internal wave solutions are constructed for a selected number of simple depth functions $d$.
\end{abstract}

\medskip




\section{Introduction}\label{sec:Introduction}



Internal gravity waves form the final chapter of a classic book on ``Waves in Fluids''~\cite{Li78}.
Equation (22) at~\cite{Li78} states that the the upward component of the mass flux, $q$ satisfies
$$\Delta (\frac{\partial^2 q}{\partial t^2}) = -N(z)^2\left(\frac{\partial^2 q}{\partial x^2}+ \frac{\partial^2 q}{\partial y^2}\right) ,
$$
where $\Delta$ is the 3-dimensional Laplacian, and $z$ is the vertical coordinate.
Here $N(z)$ is the \BV frequency.
For 2-dimensional flows, i.e. no $y$ dependence, there is a stream function, and
several problems of physical interest involve solutions of the form 
$q(x,z,t)=\sfun(x,z)\exp(i\omega t)$,
and when, additionally, the \BV frequency is constant, $\sfun$ satisfies the hyperbolic equation~(\ref{eq:pdeD}). 

The problem we treat in this paper - standing internal waves - is not well posed, and,
in particular, solutions when they exist are not unique.
The same \pde but with different boundary conditions 
describes
 two-dimensional internal waves generated by an oscillating cylinder in a uniformly stratified fluid:
and a few comments on this are given in  our \S\ref{sec:Discussion}.
A photograph of the wave pattern of this is given in Figure~76 on page~314 of~\cite{Li78} and
a diagram indicating the beams of internal waves is given in Figure~2 of~\cite{Hu97}.
The characteristic directions of the \pde are very evident.
For our standing wave problem,
once again the characteristic directions are often evident in the flow fields: see,
for example,  our Figure~\ref{fig:one} and other publications on the subject,
including photographs of experiments.

For general plane domains standing waves are treated in~\cite{AK98}. 
In this paper we specialise to fluid domains confined by a flat surface $z=0$ and a bottom boundary $z=-d(x)$  for a given non-negative depth function $d$. 
Exact solutions for certain depth functions $d$ are known, e.g. 
Wunsch's solution for a subcritical wedge \cite{Wu68}, 
Barcilon's solution in a semi-ellipse \cite{Ba68} and 
a self-similar solution in a specific trapezoid \cite{Ma09}, 
among many others. 
It is known that analytical solutions to the governing differential equation with Dirichlet boundary condition can be constructed from functions $f$ which satisfy the functional equation 
$$
f\left(x+\frac{d(x)}{\nu}\right)=f\left(x-\frac{d(x)}{\nu}\right) + Q
\eqno{({\rm FEd}(Q))}
$$
for $\nu >0$ and $Q$ given constants.  
This functional equation has been used for internal wave studies for several decades: see~\cite{ML95} and references therein. The physical interpretation of $Q$ non-zero is a constant mass-flux through the domain and it is considered in \cite{Ma11}  in the context of tidal conversion. 
The zero-flux boundary condition $Q=0$ is the physical condition appropriate to standing waves (and blinking modes) and is the main topic of this article. 
It has been noticed by \cite{MM71} (their Theorem 2) and \cite{Sa76} that there are reformulations of FEd($0$) such that one can associate solutions to FEd($0$) with solutions to FEd($Q\neq 0$). 
However, to date, very little 
use of advantages associated with these reformulations 
seems to have been made in the construction of analytical internal wave solutions. 

For a large class of depth functions $d$ one can invert the arguments in the functional equation FEd($Q$) and formulate them as the functional equation FET($Q$) presented in \S\ref{sec:FunctionalEquations}, 
which corresponds to a special case of Schr\"{o}der's functional equation for $Q=0$ and Abel's functional equation $Q\neq0$.  
Schr\"{o}der's and Abel's functional equations are well-studied functional equations~\cite{Ku68}. 
In this article known properties of these functional equations are put into context for the construction of internal waves. 
A selection of analytical internal wave solutions constructed from solutions to these functional equations is presented. 
Besides the application to internal waves,
there are other wave phenomena described by the same boundary-value problem:
we mention some of these at the end of \S\ref{sec:pde}.

The structure of this paper is as follows.
In~\S\ref{sec:pde} we present the partial differential equation boundary-value problem
that models the internal waves and
in~\S\ref{sec:FunctionalEquations} we present the corresponding functional equations.
We present in~\S\ref{sec:Wunsch} Wunsch's solution for a subcritical wedge,
and follow this in \S\ref{sec:subcritical} 
with various solutions for standing waves with everywhere subcritical bottom profiles.
Our treatment in~\S\ref{sec:partSuper} and in~\S\ref{sec:Involutions} indicates results for bottom profiles
that have some supercritical parts.
The latter of these two sections, \S\ref{sec:Involutions},
treats a particularly simple solution method appropriate when $d$ is related in a certain way to involutions.
We are confident that the methods allow for further development:
related problems where they might be used are presented in~\S\ref{sec:Discussion}.

\section{Internal wave differential equation} \label{sec:pde}

Let the bottom topography $d(x)$ be a positive function defined on the open connected interval $I=[b_-,b_+]\subset \mathbf{R}$. If $b_{\pm}$ are finite, then $d(b_{\pm})=0$. 
Define the simply-connected open domain $D$ in the plane by
$$ D= \{(x,z)\in\mathbf{R}^2\, | \, b_{-}<x<b_{+}, -d(x)< z< 0\},$$
with $x$ and $z$ representing the horizontal and vertical coordinates respectively. 
For a constant \BV frequency, the streamfunction $\sfun$ of small-amplitude internal waves in $D$ is governed by
\begin{equation}
\begin{split}
\begin{aligned}
\,& \frac{\partial^2\sfun}{\partial x^2}- \nu^2\frac{\partial^2\sfun}{\partial z^2}=0 {\rm\qquad{in}\ } D, \\
\,&\sfun(x,0)=0 \qquad {\rm\qquad{ for}\ } b_-<x<b_+ \\
\,&\sfun(x,-d(x))=Q \ {\rm\qquad{for}\ } b_-<x<b_+
\label{eq:pdeD}
\end{aligned}
\end{split}
\end{equation}
where $\nu>0$ and $Q$ are given constants. 
A derivation of (\ref{eq:pdeD}) 
can be found in many books on fluid dynamics, e.g. 
Chapter VI \S4 on Sobolev's equation in \cite{AK98}.  
See also~\cite{ML95} equations~(2.4) and (2.5)-(2.6), the latter specifically for the case $Q=0$.
For $Q$ nonzero, see~\cite{Ma11}, in particular the 
paragraph containing his equation~(2.1).

The quantity $\nu$ can be interpreted as the inclination of the characteristics
(internal wave rays or beams) relative to the horizontal.
A point $x$ on the bottom of the domain $D$ is called {\it subcritical} if
the bottom topography function $d$ satisfies 
$|d'(x)|<\nu$, where $d'$ denotes the derivative of $d$, 
and {\it supercritical} if the reverse holds. 
If all points on the bottom are subcritical (supercritical),
then the bottom profile $d$ and the domain $D$ are each refered to as being subcritical (supercritical). 

Notice that it is always possible to stretch the $z$-coordinate such that 
$\nu$ takes the value 1 in the problem with the scaled bottom topography $d(x)/\nu$. 
In the following,
unless $\nu$ is explicitly referenced, 
the parameter $\nu>0$ is assumed to be 1.

We will consider $Q\neq 0$ when it is appropriate.
This happens when all points on the bottom are subcritical
(see \S\ref{sec:subcritical}), and in some other instances
(see \S\ref{subsec:daiHyperbola}).
For bounded domains $D$ the physical interpretation has 
(sinusoidally oscillating) sources and sinks at $(b_\pm,0)$.

Various comments are appropriate.
The standing wave solutions, i.e. those with $Q=0$, harmonic in time, 
can be used to solve initial-boundary-value problems for the Sobolev equation. 
Related problems occur in other applications, for example, 
in some theoretical physics applications (e.g.~\cite{Do01}),
and other moving boundary problems for the wave equation (e.g.~\cite{Di98}).

\section{Functional equations}\label{sec:FunctionalEquations}

The functional equations in this paper are all linear;
the $Q=0$ case being homogeneous.
Some properties hold for any $Q$ zero or nonzero.
If one has a solution $f$ then $f+c$ is also a solution
for any constant $c$.
Suppose $f_0$ and $f_1$ are solutions at the same $Q$.
The minimum of $f_0$ and $f_1$ is also a solution.
The convex combination $(1-t)f_0+tf_1$ is also a solution.
Consequences of these are used without further comment in this paper.

\subsection{The `extension of $f$' to $\sfun$}\label{subsec:extension}
Assume a solution of the differential equation in (\ref{eq:pdeD}) 
is represented by 
\begin{equation}
\sfun(x,z) = f\left( x-\frac{z}{\nu}\right)- f\left(x+\frac{z}{\nu}\right) \qquad {\rm for} \ (x,z)\in D
\label{eq:dalam}
\end{equation}
for some differentiable real function $f$. 
The boundary condition $\sfun(x,d(x))=Q$ is satisfied if $f$ satisfies the functional equation FEd($Q$) given in~\S{1}. 
Note that $\sfun(x,0)=0$ is already satisfied by the definition~(\ref{eq:dalam}). 


With $\sfun$ defined from~(\ref{eq:dalam}), $\sfun$ will inherit smoothness properties from $f$.
Piecewise linear functions $f$ will produce piecewise linear $\sfun$.

We have used the term `extends'  merely to indicate the following.
Given a function $f$ defined on an interval $(c_-,c_+)$
one can view equation~(\ref{eq:dalam}) as extending the one-dimensional domain $(c_-,c_+)$
to a domain in the plane. 
(Strictly speaking $f$ itself might better be thought of
as extending to the hyperbolic conjugate of $\sfun$~\cite{ML95}
as this is such that its restriction to $z=0$ is, except for a factor of 2, the function $f$.)
This extension defines the function $\sfun$ in the triangle in $z\le{0}$ 
with its other sides the characteristics through $(c_\pm,0)$,
namely the lines $z=c_\pm\mp\nu{x}$.
When $d$ is everywhere subcritical, we can take $c_\pm=b_\pm$ and, 
when both $b_+$ and $b_-$ are bounded, the triangle so formed contains the whole of the domain $D$.
The extension via~(\ref{eq:dalam}) might well lead to a $\sfun$ defined over a larger set than the domain $D$. 
In the case $Q=0$, the curve $z=-d(x)$ is then a nodal curve of $\sfun$ defined over the larger set.

Suppose now that $b_-=-b_+$. 
When $f$ is an even function the corresponding $\sfun$ is odd in $x$.
When $f$ is an odd function the corresponding $\sfun$ is even in $x$.

\subsection{The forward map $T$}\label{subsec:forwardMap}

Define the functions $\delta_{\pm}:= x\pm d(x)/\nu $. 
If the $\delta_-$ in FEd($Q$) is invertible, then one can (provided the domain of $\delta_+$ includes the image of $\delta_-^{-1}$) define the  map $T_+:=\delta_{+}\circ \delta_{-}^{-1}$ and rewrite the functional equation FEd(Q) as the functional equation
$$
f(T_+(x))=f(x) + Q.
\eqno{({\rm FET}_+(Q))}
$$
In the same way, when appropriate conditions are satisfied, defining  the  map $T_-:=\delta_{-}\circ \delta_{+}^{-1}$,
one is led to the functional equation $f(x)=f(T_-(x)) + Q$.
Let $d(b_{\pm})=0$ for the remainder of this section, so that
$\delta_\pm(b_\pm)= b_\pm$.
The domains of both $\delta_-$ and $\delta_+$ are the same as the domain of $d$ namely $[b_-,b_+]$,
It remains to specify the domains of $T_+$, $T_-$ and of $f$.
It is simplest to consider a subcritical bottom $d$.
Then 
(i) both $\delta_-$ and $\delta_+$ are monotonic increasing so invertible,
(ii) the maps $T_\pm$ are bijective on $[b_-,b_+]$ -- in fact increasing on $(b_-,b_+)$ with $T_\pm(b_\pm)=b_\pm$.
To simplify notation, where this is appropriate, we omit the subscript $+$, and
the  equation we study is
$$
f(T(x))=f(x) + Q.
\eqno{({\rm FET}(Q))}
$$
For more on the case of subcritical bottoms, see the beginning of \S\ref{sec:subcritical}. 
Partly or entirely supercritical domains are more complicated: see \S\ref{sec:partSuper}. 

There are geometric and physical relations between the functions $d$ and $T$. 
A rightwards ray starting from $(x,0)$ reflects from a subcritical bottom $d$ and is next incident at the top
at $(T(x),0)$. 
(For partly supercritical bottoms, we view $(T(x),0)$ as the point where the reflected ray -- possibly prolonged through the bottom profile -- meets $z=0$, possibly with $T(x)>b_+$.)
The reflection at the bottom takes place halfway between $x$ and $T(x)$ along the $x$-coordinate and at the depth $-\nu \frac{T(x)-x}{2}$, so
\begin{equation}
d\left(\frac{x+T(x)}{2}\right) =\nu \frac{T(x)-x}{2} .
\label{eq:Tfromd}
\end{equation}
From this, with
$$ X=\frac{x+T(x)}{2}, \qquad
T(X-\frac{d(X)}{\nu})= X+\frac{d(X)}{\nu} .
$$
Provided the range of $T$ is a subset of the domain of $T$,
repeated composition -- iterates of $T$-- can be defined.
When $T$ is (strictly) increasing, with $T(b_+)=b_+$,
repeated compositions of the map $T$ applied to any $x\in(b_-,b_+)$ give a sequence $\{T^{[k]}(x)\}_{k\in \mathbf{N}}$ which converges to the fixed point $T(b_+)=b_+$ for $k\rightarrow \infty$. 
Similarly, when $T(b_-)=b_-$, one gets a sequence $\{T^{[-k]}(x)\}_{k\in \mathbf{N}}$ converging to $T(b_-)=b_-$ for repeated compositions of the inverse map $T^{[-1]}$ to any $x\in(b_-,b_+)$. 

\subsection{Schr\"{o}der functional equation {\rm FET(}$0${\rm)}}\label{subsec:schroder}

Equation~FET($0$) is a special case of the Schr\"{o}der functional equation 
$$f(T(x))=s\cdot f(x)$$ for $s=1$~\cite{Ku68}. 
This subsection presents a few properties of solutions to~FET($0$). 
A comprehensive list of known properties of Schr\"{o}der functional equation - sometimes also referred to as Schr\"{o}der-Konig's functional equation - can be found in Chapter VI of \cite{Ku68}.

One comment on the case $s>0$ is appropriate (and will be used in \S~\ref{subsec:Sandstrom}:
see equation~(\ref{eq:SaAbel})).
The following old result is standard: see, for example, \cite{Ku68} p163.
\begin{theorem}
\label{thm:S1} 
If $f$ is a positive solution of
the Schr\"{o}der functional equation $f(T(y))=s\cdot f(y)$ for $s>0$, $s\ne{1}$,
then $a(x)=\log(f(x))/\log(s)$ is a solution of the Abel equation {\rm FET(}$1${\rm )}.
\end{theorem}

Some properties of solutions of FET($0$) are easy to see.
If $T$ is not the identity function $T(x)=x$ (or equivalently if
$d$ is not the zero function), no solution of
FET($0$) (or of FEd($0$)) can be monotonic.
Hence any solution must have a local maximum or minimum
in $(b_-,b_+)$.
The solutions we present for $f$ have various numbers of
maxima and minima -- sometimes finitely many, e.g.
\S\ref{subsec:Barcilon},
sometimes countably infinitely many, e.g. the domains treated in
\S\ref{sec:subcritical}.

\begin{theorem}{
\label{thm:S2} 
If $f:I\rightarrow f(I)\subset \mathbf{R}$ is a solution to {\rm FET(}$0${\rm)} and $F$ is any real function whose domain contains the image $f(I)$ of $f$, then the composition $F\circ f$ is also a solution to {\rm FET(}$0${\rm )}.}
\end{theorem}
\par\noindent{\it Proof.}\ If $f$ is a solution of FET($0$), then $f(x)=f(T(x))$. $F$ works on the image $f(I)$ of $f$, so it follows directly that $F(f(x))=F(f(T(x))$. This shows that the composition $F\circ{f}$ also satisfies FET($0$) and completes the proof. \\

The nodal curves for $\sfun_f$ associated with $f$ according to (\ref{eq:dalam}) remain nodal curves for $\sfun_{F\circ{f}}$ associated with $F\circ{f}$. There may be more nodal curves for $\sfun_{F\circ{f}}$ unless $F$ is invertible. \\

So if a solution to Schr\"{o}der's functional equation FET($0$) exists, then it is not unique - and one can be more constructive on this point: one is free to choose a function on some subset $I_0$ of the interval $I$ on which FET($0$) must hold. 
This subset $I_0$ is refered to as a {\it fundamental interval}~\cite{ML95}. 
Once a choice for a solution $f$ on some fundamental interval $I_0$ is made, then $f$ is uniquely defined on all of $I$. 
Notice that a solution $f$ to FET($0$) takes the same value for each element of the set $\{T^{[k]}(x)\}_{k\in \mathbf{Z}}$ for each $x\in(b_-,b_+)$. 
So if $f(x)$ is prescribed for one $x\in\{T^{[k]}(x)\}_{k\in \mathbf{Z}}$, then so it is for the entire set $\{T^{[k]}(x)\}_{k\in \mathbf{Z}}$. Together with the property $T(x)>x$ it shows that $I_0=[x_0, T(x_0) )$ is a fundamental interval for any $x_0\in(b_-,b_+)$. Such a connected fundamental interval (with $x_0=0$) is considered at the beginning of \S\ref{sec:subcritical}. Be aware that 
it is not necessary for a fundamental interval $I_0$ to be a connected. 

The solvability of Schr\"{o}der functional equations FET($0$) depends crucially on the property $T^{[k]}(x)\neq x$ for all $x$ in the open interval on which FET($0$) holds and for every positive $k\in\mathbf{N}$ \cite{Ku68}. 
The following theorem deals with the consequences of fixed points of the map $T$ on the solvability of FET($0$).\\

\begin{theorem}
\label{thm:S3}{
Let $T$ be a strictly increasing continuous function on $(b_-,b_+)$ for which $T(b_{\pm})=b_{\pm}$.
Suppose also that $T^{[k]}(x)\rightarrow{b_{\pm}}$ as $k\rightarrow\pm \infty$ for $b_-<x<b_+$.
Then the only solutions of} \rm{ FET($0$)} {\it which are continuous on
the closed interval $[b_-,b_+]$ are the constant solutions.
}
\end{theorem}

\par\noindent{\it Proof.}\ Let $x\in(b_-,b_+)$. Then $f(x)=f(T(x))=\ldots=f(T^{[k]}(x))$ and $f(x)=f(T^{[-1]}(x))=\ldots=f(T^{[-k]}(x))$ for all $k\in\mathbf{N}$. By assumption on $T$ we have $T^{[k]}(x)\rightarrow{b_{\pm}}$ as $k\rightarrow\pm \infty$. Continuity of $f$ then requires $f(x)=f(T^{[k]}(x))\rightarrow f(b_{\pm})$ for $k\rightarrow \pm \infty$. So for $f$ to be well-defined for all $x\in[b_-,b_+]$ it is required that $f(x)=f(b_{+})=f(b_-)$, which allows constant solutions only. \\

\subsection{Abel functional equation {\rm FET(}$Q\neq0${\rm )}}\label{subsec:Abel}
Abel's functional equation corresponds to FET($Q$) for $Q\neq 0$. 
In some theoretical physics papers, e.g.~\cite{Do01}, it is called Moore's equation.
The physical interpretation of $Q\neq 0$ is a constant non-zero flux $Q$ through the bottom $z=-d(x)$. 
Mathematically one can treat $Q$ as a non-zero constant and associate it with the no-flux condition $Q=0$ of Schr\"{o}der's functional equations FET($0$), as motivated in the following observation.

Any solution $f$ to the Schr\"{o}der's functional equation FET($0$) has to be identical on the endpoints $x_0$ and $T(x_0)$ of a connected fundamental interval $I_0=[x_0,T(x_0))$. 
This is the motivation to consider any solution $f$ to FET($0$) to be a composition of a periodic function $P$ with an argument function $a$. 
The function $f(x)=P(a(x))$ with $P$ having period $Q>0$ then satisfies the FET($0$) if and only if the argument function $a$ satisfies one of the functional equations
\begin{equation}
a(T(x))=a(x)+Q\cdot n {\rm\qquad{for}\ n\in\mathbf{Z} }.
\label{eq:AFE}
\end{equation}
It is always possible to scale $a(x)$ such that $Q=1$. 

The fundamental interval introduced in the previous subsection applies in the same way to Abel's functional equation, 
e.g. if a solution exists, then it is uniquely determined if and only if it is prescribed on a fundamental interval. 
(See the beginning of \S\ref{sec:subcritical} for an existence result.)

\begin{theorem}
{\label{thm:A1}
Let $a \in C^1$ be a strictly increasing solution of~FET(1).\\
(1) The general solution $a_{gen}$ of~FET(1) is given by 
$$a_{\rm gen}(x)= a(x) + P(a(x)) $$
where $P$ is a periodic function with period 1.\\
(2) If $a^*$ is another strictly increasing $C^1$ solution of~FET(1) then there exists some periodic function $P$ with period 1 such that $P'(x)>-1$ for all $x$ and
\begin{equation}
a^*(x)= a(x) + P(a(x)).
\label{eq:szAbel}
\end{equation}
Conversely any $a^*$ of the form~(\ref{eq:szAbel})
is an invertible solution of FET(1).
}
\end{theorem}

Part (1) is Theorem 1 of~\cite{Ne98}. Part (2) is from
\cite{Sz98} 
who attributes it to Abel (1881). 
Part (2), with its condition $P'(x)>-1$ is developed for $C^k$ solutions in
Theorem 2 of~\cite{Ne98}, with further development in his Theorem 3.\\

\begin{theorem}
\label{thm:A2}{
Let $a$ and $f$ be $C^1$ solutions to respectively {\rm FET($1$)} and {\rm FET($0$)} on $I$. Assume further that $a$ is injective and $T:I\rightarrow I$ bijective. Then there exists some periodic function $P$, with period 1, such that $f(x)=P(a(x)).$}
\end{theorem}
\par\noindent{\it Proof.}\ For $a$ as given in the theorem there exists an inverse $a^{-1}$ on the image $a(I)$ of $a$. Define $P=f\circ a^{-1}: a(I)\rightarrow \mathbf{R}$. 
It is easy to verify that this function $P$ satisfies $f=P\circ a$. The claim is that this function is periodic. The functions $f$ and $a$ satisfy $f(T(x))=f(x)$ and $a(T(x))=a(x)+1$. This gives $P(a(T(x)))=f(T(x))=f(x)=P(a(x))=P(a(T(x))-1)$. The assumption that $T$ maps its domain $I$ bijectively onto itself, $T(I)=I$, gives that $P(x)=P(x-1)$ for all $x\in a(I)$, which shows that $P$ is periodic with period $1$. \\

A direct consequence of Theorem~\ref{thm:A2} is that for subcritical bottom topographies all continuous solutions to FET($0$) are constructed by applying the set of all continuous periodic functions with period $Q>0$ to any continuous injective solution to FET($Q$).
\begin{theorem}
\label{thm:A3}{
Given a strictly increasing continuous map $T$ on $(b_-,b_+)$ with $T(b_{\pm})=b_{\pm}$, some fundamental interval $I_0=[x_0, T(x_0))$ and a strictly increasing continuous function $a_0$ on $I_0$, then the unique continuous solution $a$ to FET($Q$) with $a=a_0$ on $I_0$ and $Q=a_0(T(x_0))-a_0(x_0)$  satisfies

\begin{equation}
a(x)=a_0(T^{[-k]}(x))+k Q
\label{abelalgorithm}
\end{equation}
for all $x\in I_k:=[T^{[k-1]}(x_0),T^{[k]}(x_0))$ and $k\in\mathbf{Z}$.
}
\end{theorem}
This theorem is a special case of Theorem 4.1 in \cite{La07}, which proves that $a(x)=a_0(T^{[-k]}(x))+k Q$ for $x \in I_k$ if $a$ is continuous solution satisfying FET($Q$). In \cite{La07} the function $a_0$ satisfying $a_0(T(x_0))-a_0(x_0)=Q$ is assumed to be linear, which is in fact not necessary for the proof. \\
The solution $a(x)$ to FET($Q$) is clearly continous in all points $x$ in the interior of some interval $I_k$. For the boundary points $x_k:=T^{[k]}(x_0)$ study the limits $x\rightarrow x_k$ for $x>x_k$ and $x<x_k$:
If $x>x_k$, $x\in I_k=[x_k,x_{k+1})$, then
$$\lim\limits_{x \to x_k} a(x)=a_0(T^{[-k]}(x_k))+k Q=a_0(x_0)+kQ.$$
For $x<x_k$, $x\in I_{k-1}=[x_{k-1},x_k)$ it follows that
$$\lim\limits_{x \to x_k}a(x)=a_0(T^{[-k+1)]}(x_k))+(k-1) Q=a_0(T(x_0))+(k-1)Q.$$ 
These two expressions are equal because $a(T(x_0))=a(x_0)+Q$ by the definition of $Q$. \\
To prove uniqueness observe that for every $x\in (b_-,b_+)$ there exists a unique $k\in\mathbf{Z}$ such that $x\in I_k$ because $\bigcup \limits_{k=-\infty}^{\infty} I_k=[b_-,b_+]$ and all $I_k$ are disjunct. So for every $x \in (b_-,b_+)$ the function $a(x)$ is uniquely defined by the expression (\ref{abelalgorithm}) since $T$, $a_0$ and $Q$ are given.

\subsection{Comments on FEd($0$) and FEd($Q$)}\label{subsec:FEd}

Some of the results of
\S\ref{subsec:schroder} and of \S\ref{subsec:Abel}
have analogues for the equations 
FEd($0$) and FEd($Q$) respectively.
In particular, we remark that if $P$ is a periodic function with period $Q>0$ and $a$ solves FEd($Q$), then the composition
$P_Q\circ{a}$ solves the Schr\"{o}der-like equation FEd($0$).

\section{Wunsch's solution: subcritical wedge\label{sec:Wunsch}}

Let $b_-=-\infty$, $b_+\in \mathbf{R}$ and $\nu=1$. 
For a subcritical wedge $d(x)=\tau(b_+-x)$ with $\tau\in(0,\nu)$ 
the map $T$ is the linear function 
$T(x)=px +s$ where $p=\frac{1-\tau}{1+\tau}$ and $s= b_+\frac{2\tau}{1+\tau}$. The Schr\"{o}der functional equation FET($0$)
\begin{equation}
f(p x+s)=f(x) {\rm\qquad{for}\ }x< b_+ 
\label{eq:Wunsch1}
\end{equation}
can be formulated as the Abel's functional equation FET($1$) under the assumption $f=P\circ a$ with $P$ any period-$1$ function:
\begin{equation}
a(p x+s)=a(x)+1 {\rm\qquad{for}\ } x< b_+ .
\label{eq:Wunsch2}
\end{equation}
A continuous, stricktly increasing 
solution to (\ref{eq:Wunsch2}) is given by 
$a(x)={\log(-x+b_+)}/{\log(p)}$.
So the Schr\"{o}der functional equation (\ref{eq:Wunsch1}) is solved by 
functions 
$$f(x)=P\left(\frac{\log(-x+b_+)}{\log(p)}\right)$$
for any arbitrary continuous period-1 function $P$. 

The solution given by \cite{Wu68} had $P$ as a sine or cosine function.
The nodal curves which intersect $z=0$ in these solutions 
are hyperbolae. Of course there are many other periodic functions.
For certain piecewise exponential $P$ all the nodal curves are straight lines:
for appropriate $P$ some nodal lines are vertical straight lines.
This makes a connection with 
this section and~\S\ref{subsec:isosctriangle}.

\section{Symmetric domains with subcritical bottom profiles}\label{sec:subcritical}


Our treatment of the functional equations in~\S\ref{sec:FunctionalEquations} deliberately avoided
general existence matters as these can be rather intricate, except in the context of subcritical bottoms.
The existence result in the next paragraph is stated as it
 provides a lead-in to~\S\ref{subsec:isosctriangle}.

In the existence result below we have a genuine interval as
a fundamental interval.
(That this is not always the case is mentioned 
in~\S\ref{subsec:schroder}.)
For a symmetric domain, take as the domain of $x$ the interval $[b_{-},b_{+}]=[-b,b]$ for some $b>0$.
The following is stated in~\cite{Ne98} (giving references for the proof, including~\cite{Ku68}).\\
\begin{theorem}{ {\label{thm:subcrit1}}
$T$ is a continuous strictly increasing real-valued function defined on a half- open interval $[0,b)$, $0<b\le\infty$,\\
$T([0,b))=[c,b)$ with $c>0$, (so we can extend, by continuity, the domain of $T$ so $T(b)=b$) and\\
$T(x)>x$ for $0\le{x}<b$\\
then there exists a solution for {\rm FET(1)}. Furthermore under the above conditions,
there is a unique solution $a$ with prescribed values on the interval $[0,T(0))$.
If, moreover, it is continuous on $[0,T(0))$ and (taking the limit from above)
$$ \lim_{x\rightarrow{T(0)}} a(x)  
=a(0)+1 $$
then $a$ is continuous on $[0,b)$.
}
\end{theorem}
All the conditions on $T$ above are satisfied by the forward maps $T$ of symmetric domains with subcritical bottom profiles.
{ (A hydrodynamic interpretation is that, for a given bottom profile $d$,
there is a solution for all $\nu$ satisfying $\nu>{\rm max}(|d'(x)|)$.)}\\

Any such solution $a$ necessarily tends to minus infinity as $x$ tends to $b_-$,
and to plus infinity as $x$ tends to $b_+$.
(If $a$ were to be continuous on the closed interval $[b_-,b_+]$ the solutions of the
Schr\"{o}der equation generated from it could also be continuous, contradicting Theorem~\ref{thm:S3}.)

In the context of the symmetric domains and $Q\ne{0}$ our main interest is in odd solutions $a$.

\subsection{Subcritical isosceles triangle}\label{subsec:isosctriangle}

In this section construct  all possible solutions to FET($0$) for the isosceles triangle with bottom topography function $d(x)=\tau(1-|x|)$ with $\tau\in(0,1)$ for $x\in(b_-,b_+)=(-1,1)$ and $\nu=1$ are constructed. According to Theorem~\ref{thm:A2} one can construct all solutions $f$ to FET($0$) via the relation $f=P\circ a$ with $P$ any periodic function with period $Q$ (=length of $I_0$) and $a$ a continuous, strictly increasing solution to Abel's functional equation FET($Q$). The goal is therefore to construct one solution to FET($Q$) for some $Q\neq 0$ using the expression (\ref{abelalgorithm}) from Theorem (\ref{thm:A3}). 
The map $T=\delta_+\circ \delta_-^{-1}$ and its inverse $T^{[-1]}$ associated with $\delta_{\pm}=x\pm d(x)$ are given by
\begin{equation}
\begin{split}
\begin{aligned}
\,&T(x)=p^{-1} x+s_- \,&{\rm\qquad{for}\ -1\leq x\leq -\tau}\\
\,&T(x)=p x+s_+ \,&{\rm\qquad{for}\ -\tau\leq x\leq +1}\\
\,&T^{[-1]}(x)=p x-s_+ \,&{\rm\qquad{for}\ -1\leq x\leq +\tau}\\
\,&T^{[-1]}(x)=p^{-1}x-s_- \,&{\rm\qquad{for}\ +\tau\leq x\leq +1}
\label{eq:T}
\end{aligned}
\end{split}
\end{equation}
where $p=\frac{1-\tau}{1+\tau}<1$, $s_+=\frac{2\tau}{1+\tau}$ and $s_-=\frac{2 \tau}{1-\tau}$. 
A fundamental interval is given by $I_0=[-\tau,\tau)$, as can be verified by checking that $T(-\tau)=\tau$. 
Repeated compositions of function $T$ and its inverse $T^{[-1]}$ map this fundamental interval $I_0$ onto the intervals $I_k:=T^{[k]}(I_0)$, $k\in\mathbf{Z}$. 
So for  $x\in I_k$ and $k\leq -1$ a solution $a(x)$ to the Abel equation $FET(Q)$ is given by $a(x)=a_0(T^{[k]}(x))-k Q$ where $a_0$ is an arbitrary strictly increasing choice for $a$ on $I_0$ which satisfies $a_0(\tau)-a_0(\tau)=Q$. Similarly for $k\geq 1$ and $x\in I_k$ one gets $a(x)=a_0(T^{[-k]}(x))+kQ$.\\
Compositions of the maps $T$ and $T^{[-1]}$ give 
\begin{equation}
\begin{split}
\begin{aligned}
\,&T^{[k]}(x)=1+p^{-k}(x-1) \ \,&{\mbox{ for }} \ -\tau<x \\
\,&T^{[-k]}(x)=-1+p^{-k}(x+1) \ \,&{\mbox{ for }}  x<+\tau . 
\label{eq:Tk}
\end{aligned}
\end{split}
\end{equation}
For the simple choice $a_0(x)=x$ on the fundamental interval $I_0$, which implies $Q=a_0(\tau)-a_0(-\tau)=2\tau$, the continuous solution $a$ is given by
\begin{equation}
\begin{split}
\begin{aligned}
\,&a(x)=p^{-n}(x-1)+1+2\tau n {\rm\qquad{for}\ x  \in I_n, n\in\mathbf{N}}\\
\,&a(x)=p^{-n}(x+1)-1-2\tau n {\rm\qquad{for}\ x  \in I_{-n}, n\in\mathbf{N}}.
\label{eq:sol}
\end{aligned}
\end{split}
\end{equation}



\begin{figure}
\centerline{\includegraphics[height=7cm,width=13cm]{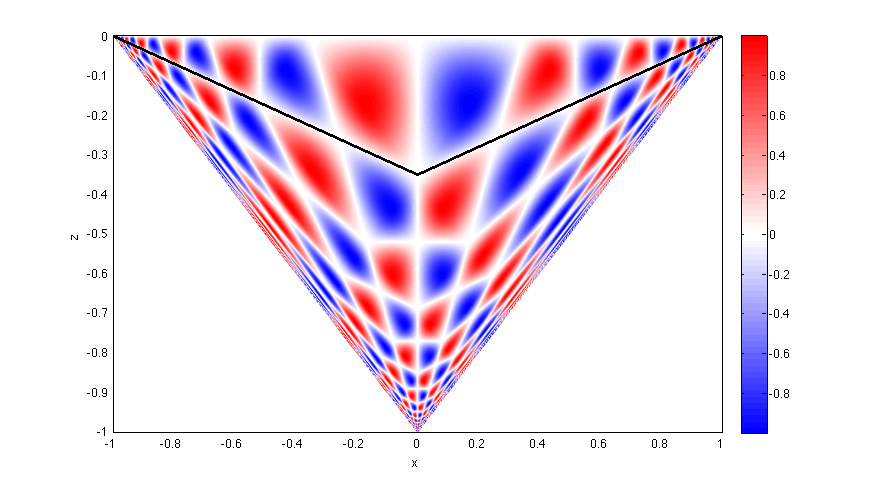}}
\caption{This figure shows the analytical streamfunction solution for $\tau=0.35$ with $P(x)=\cos(\frac{\pi}{\tau} x)$. The bottom of the isosceles triangle is indicated by the black line. All streamfunction values $z<|x|-1$ are set to zero.}
\label{fig:tri}
\end{figure}

In Figure \ref{fig:tri} a continuously differentiable streamfunction solution $\Psi(x,z)=f(x-z)-f(x+z)$ for the choice $P(x)=\cos(\frac{\pi}{\tau} x)$ is presented. The black line shows the bottom $d(x)=\tau(|x|-1)$. It appears that there are many nodal curves. 
The plotted solution is also a solution for many bottom topographies, including partly and entirely supercritical bottom topographies. It is speculated that some of these nodal curves are independent of the choice of the periodic function $P$, e.g. streamfunction solutions to the bottom topographies along these isoclines can be constructed from $f=P\circ a$ for arbitrary period-$2\tau$ function $P$ and $a$ satisfying (\ref{eq:sol}). 

\subsection{Subcritical symmetric hyperbolae}\label{subsec:Sandstrom}

\subsubsection{Symmetric hyperbolic lens}
\begin{figure}
\centering
$$
\includegraphics[height=7cm,width=13cm]{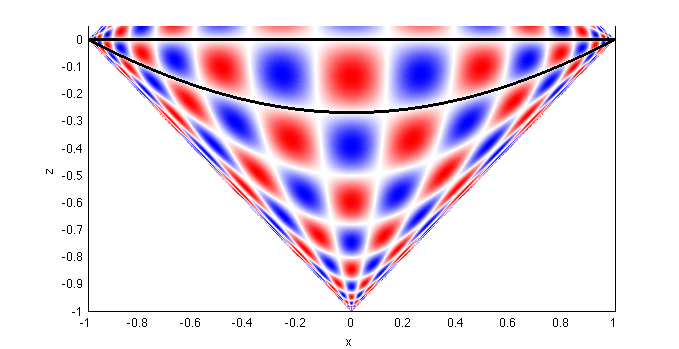}
$$
\caption{The streamfunction soltuion $\Psi(x,z)=f(x-z)+f(x+z)$ is plotted with $f$ being the composition of $P(x)=\sin(\frac{2\pi}{\rm{arctanh}(1/c)}x )$ for $c=2$ and $a(x)=$arctanh$(x)$ (which solves FET($a(1/c)$)). The color bar is as in Figure \ref{fig:tri}.}
\label{fig:hypSymSubcrit}
\end{figure}

Again, set $\nu=1$.
For the subcritical bottom topography 
\begin{equation}
d(x) = c -\sqrt{c^2-1+x^2} \qquad {\rm for\ } -1<x<1 \qquad {\rm with\ } c>1
\label{eq:Sa76}
\end{equation}
the corresponding map $T$ is given by
\begin{equation}
\begin{split}
T(x)=\frac{1+ c x}{c+x} 
= x + \frac{1-x^2}{c+x} \qquad {\rm for\ } -1<x<1 .
\end{split}
\end{equation}
The map $T$ is fractional linear.
Defining another fractional linear map $r$ and motivated by the fact that compositions of fractional linear maps are fractional linear,
$$r(x)= \frac{1+x}{1-x} {\mbox{\rm\ \ \ gives\ \ }}
r(T(x))= r\left(\frac{1}{c}\right) r(x) . $$
(The function $r$ satisfies a Schr\"{o}der functional equation with
$s=r(1/c)$ positive.)
Take logs of $r(x)$ and notice that $a(x)=\frac{1}{2} \log(r(x))=$arctanh$(x)$ satisfies
\begin{equation}
a(T(x)) = a(x) + a(\frac{1}{c}).
\label{eq:SaAbel}
\end{equation}
This solution has been suggested by~\cite{Sa76}. The solution $a(x)=$arctanh$(x)$ is injective on the fundamental interval $I_0=[0,\frac{1}{c})$ because $\frac{1}{c}<1$. So according to Theorem~\ref{thm:A2} all solutions $f$ to FET($0$) can be derived by applying arbitrary periodic function $P$ with period 
$a(\frac{1}{c})=\frac{1}{2}\log\left(\frac{1+c}{-1+c}\right)$ to $a(x)$: $f(x)=P($arctanh$(x))$. The streamfunction solution for a sinusoidal choice for $P$ is shown in Figure \ref{fig:hypSymSubcrit}. 

There are infinitely many nodal curves intersecting $z=0$ at points in $-1<x<1$.
Modes with different numbers of cells stacked vertically are easily constructed.

\subsection{Some other subcritical bottom profiles}\label{subsec:othersubcrit}

The entries in
the table
indicate some other subcritical bottom profiles for 
which we have solutions (with $\nu=1$).
The column headed $a$ gives solutions of the Abel functional
equation for the given $T$ (from which one can generate all
standing-wave solutions).
A banal comment -- useful when both $a$ and its inverse $a^{-1}$
have simple forms -- is the simple formula for $T$ given $a$
solving FET($Q$):
$$ T(x,Q) = a^{-1}\left( a(x) + Q\right), \qquad
T^{[k]}(x,Q) = a^{-1}\left( a(x) + k Q\right)\qquad
  .
$$

\noindent{}
\begin{tabular}{|| c | c | c | l ||}
\hline
$[b_-,b_+]$& $T$ & $a$ & 
Comments \\
\hline
$[0,1/2]$& $2x(1-x)$& $ 
\frac{\log \left(\frac{\log \left(1-2x\right)}{\log
   (1-2 c)}\right)}{\log (2)}$&  Unsymmetrical parabolic segment\\
 & & & \\
$(-\infty,\infty)$& See below &${\rm arcsinh}(x)$&
Symmetric hyperbolic hump\\
 & & & \\
See below& $ \frac{x}{1+x}$& $\frac{1}{x}$&
Source where a hyperbolic \\
 &  &  &slope intersects $z=0$\\
\hline
\end{tabular}

\par\noindent
$\bullet$
For the symmetric hyperbolic hump, for an appropriate value of 
$\tau$ with $0<\tau<1$,
$$T_{\tau}(x)= \frac{ (1+\tau^2) x+ 2 t\sqrt{1+x^2}}{1-\tau^2} ,\qquad
d_{\tau}(x)= \tau\sqrt{\frac{1}{1-\tau^2}+x^2} .$$
$\bullet$
The entry in the table corresponding to $a(x)=1/x$ can be viewed as a singular flow
corresponding to a dipole located at the origin.
(The domain of $a$ is no longer an interval.)
All streamlines are hyperbolas passing through the origin and located in the
wedge shapes containing $z=0$ and bounded by
characteristics through the origin.

There are many other solutions in the literature
e.g. in~\cite{Da11,Do01}.
A symmetrically placed fully submerged subcritical (isosceles) wedge will yield to the methods of~\S\ref{subsec:isosctriangle}.

\section{Some domains where part or all of the bottom is supercritical}\label{sec:partSuper}

Here we are concerned with solutions of  equation~FEd($0$)
$$f(x+ \frac{d(x)}{\nu})= f( x- \frac{d(x)}{\nu})) $$
where the function $f$ may need to be defined on a larger interval than is the function $d$.
$[b_-,b_+]\times\{-1,+1\}$: 
I.e. we are treating the case $Q=0$.
However in \S\ref{subsec:daiHyperbola},
we solve FEd($Q$) with $Q>0$ as part of the
metnod of solving FEd(0).
In this section we use the FEd formulations and in \S\ref{sec:Involutions} the FET version.
When the domain of $f$ is larger than
that of $d$ it restricts us to functions which extend to
a $\sfun$ with a domain larger than $D$ and vanishing on $z=0$ over more than
that part which is on the boundary of $D$: we may be finding just some of the
solutions of the differential equation problem~(\ref{eq:pdeD}).
By treating the problem in the form~FEd($Q$)
rather than FET($Q$) we avoid some of the difficulties associated with the lack of invertibility of one or other of 
$\delta_+$ or $\delta_-$.

There are other methods of solving the problem,
some of which are mentioned at the end of this 
section.

\subsection{Barcilon's solutions for the semi-ellipse}\label{subsec:Barcilon}

Let the bottom topography be a semi-ellipse: $d(x)=\sqrt{1-x^2}$ for $x\in(-1,1)$. 
The functional equation ~FEd($0$) then becomes
$$ 
f( x-\frac{\sqrt{1-x^2}}{\nu})- f( x+\frac{\sqrt{1-x^2}}{\nu})=0 .
$$
With this restriction the preceding functional equation can be re-written
\begin{equation}
f(\cos(\theta)-\sin(\theta)/\nu)- f( \cos(\theta)+\sin(\theta)/\nu)=0
\label{eq:circtheta}
\end{equation}
A family of solutions, involving Chebyshev polynomials is given in~\cite{Ba68}.
These solutions have been rediscovered several times, e.g.~\cite{ML95}.

\subsubsection{Reduction to a constant coefficient functional equation}

\begin{figure}
\centering
$$
\includegraphics[width=\textwidth]{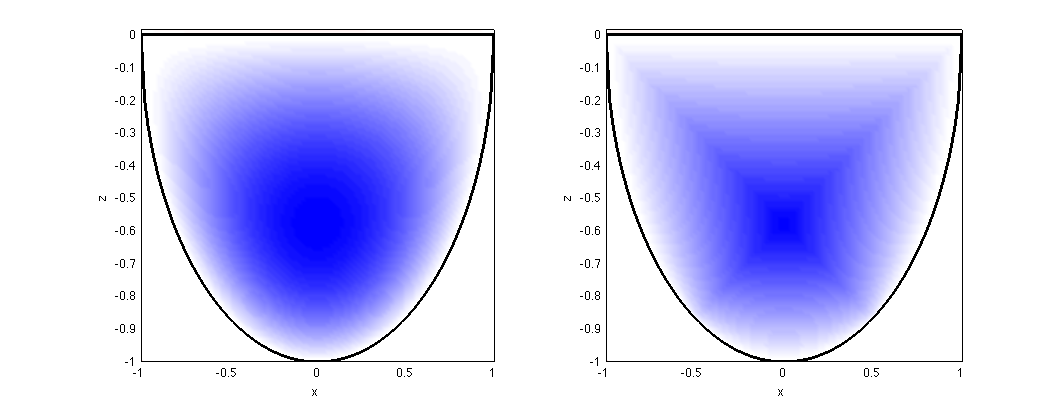}
$$
\caption{Different solutions with $k=3$, $m=1$.
At left the periodic function is $\cos$.
At right, the periodic function replaces $\cos$ with a $2\pi$ periodic even triangle wave. 
In the same way as a triangle wave can be expressed as a Fourier cosine series, the solution
at right can be represented as an infinite series superposition of polynomial solutions.}\label{fig:one}
\end{figure}

We now indicate one method to solve 
the functional equation~(\ref{eq:circtheta}), and find, 
amongst others, the Chebyshev function solutions.
We begin with seeking solutions to
$$f_+ =f( \cos(\theta)-\sin(\theta)/\nu)=f( \cos(\theta)+\sin(\theta)/\nu)=f_- . $$
Next define $\cos(\theta_\nu)=\nu/\sqrt{1+\nu^2}$.
Define also ${\tilde{f}}({\tilde{\theta}})=f(\sqrt{1+\nu^2}\cos({\tilde{\theta}})/\nu)$.
The functional equation in terms of $\tilde{f}$ is:
$${\tilde{f}}(\theta+\theta_\nu)={\tilde{f}}(\theta-\theta_\nu)$$
or, equivalently
$${\tilde{f}}(\theta)={\tilde{f}}(\theta+2\theta_\nu) .$$
This is solved, for $\tilde{f}$, by {\it any} periodic function $P$ with
period $2\theta_\nu$. 
However restrictions on $\nu$ may be
required to ensure that the extension of $f$ to $\sfun$ leads to a physically acceptable $\sfun$. 
Barcilon's Chebyshev solutions are,
with integer $m$ and $k$, from
$${\tilde{f}}({\tilde{\theta}})=\cos(\frac{m\pi{\tilde{\theta}}}{\theta_\nu})
\qquad{\rm with}\qquad
\theta_\nu = \frac{m\pi}{k} .$$
Returning to the general $2\theta_\nu$-periodic
$\tilde{f}$, having found $\tilde{f}$
we can determine $f$ as follows. Set
$$ X
=\frac{\sqrt{1+\nu^2}}{\nu}\cos({\tilde\theta})
=\frac{\cos({\tilde{\theta}})}{\cos(\theta_\nu)}, \qquad
{\tilde{\theta}}=\arccos(X\cos(\theta_\nu)) ,
$$
$$
f(X)={\tilde{f}}(\arccos(X\cos(\theta_\nu)).
$$
For Barcilon's solutions this is 
$$ \nu=\cot(\frac{m\pi}{k})\qquad
f(X)=\cos(k\arccos(\cos(\frac{m\pi}{k})X) ) $$
A couple of solutions for the lowest mode -- no interior nodal curves --
(and $\nu=1/\sqrt{3}$) are shown in Figure~\ref{fig:one}.
\medskip

For plots of some other modes, see~\cite{Ba68,ML95}.

\subsubsection{Taylor series methods for {\rm FEd($0$)} and {\rm FET($0$)} }

There are other methods that can be used to solve FEd(0) with $d(x)=\sqrt{1-x^2}$.
One can form a Taylor series about $x=0$ of each of $f(x \pm d(x)/\nu)$.
If one is to seek a polynomial solution the Taylor series is a finite sum, and furthermore
only even powers of $d(x)$ enter the equation to be solved.
It is easy to recover Barcilon's Chebyshev polynomial solutions from this approach.
One can also find other $d(x)$ which lead to polynomial $f$.
The method can also be adapted to shapes
other than the semiellipse, finding rational functions $f$, and to
solving the Abel's functional equation ($Q$ non-zero) not merely the $Q=0$ Schr\"{o}der functional equations.

\subsubsection{A forward map $T$ with range bigger than $[-1,1]$}

$T$ 
(determined using equation~(\ref{eq:Tfromd}))
is 
$$T(X)= \frac{2 \sqrt{1-\nu ^2 \left(X^2-1\right)}+\left(\nu ^2-1\right) X}{\nu ^2+1} $$


Barcilon's Chebyshev solutions of $f$ satisfying $f(X)=f(T(X))$ are readily verified.
(An easy example is $f(X)=2X^4-4X^2+1=T_4(X\sin(\pi/4))$ corresponding to $\nu=1$ and $T(X)=\sqrt{2-X^2}$.
Here $T_4$ denotes the Chebyshev polynomial of degree 4.)

\subsection{Dai's solutions for hyperbolae}
\label{subsec:daiHyperbola}

The case of a hyperbolic bottom profile $d(x)=r/x$ for $x>0$
is treated in~\cite{Da08}. 
One readily verifies that
FEd($Q$)
$$ a\left(x+\frac{r}{\nu x}\right) = a\left(x-\frac{r}{\nu x}\right)+ Q\ \
{\mbox{\rm is solved by }}\
a(x)= \frac{Q \nu x^2}{4 r}. $$
The streamfunction associated with this $a$
has fluid entering from $(\infty,0)$ and
exiting via $(0,-\infty)$.

%
%
%

In this case it happens that the problem can be recast using
the forward map $T(x)=\sqrt{4 r/\nu + x^2}$ for $x>0$ into an
Abel equation FET(Q).
The solution appears elsewhere. For example, \cite{Do01}, near his equation (9), 
gives the solution with
$$d(x)=\frac{1}{d_0+r x}
{\mbox{\ {\rm and}\ }} \nu=1, \qquad a(x)= -\frac{Q}{2} (d_0 x +\frac{r}{2} x^2 ) . $$

\begin{figure}
\centering
$$
\includegraphics[height=7cm,width=13cm]{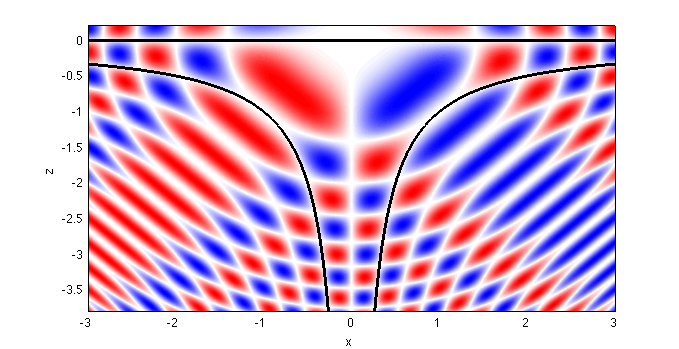}
$$
\caption{Dai's streamfunction solution for hyperbolic bottom profile $d(x)=1/|x|$ corresponding to the solution $f(x)=\cos(\frac{\pi}{2} x^2)$ to FET($0$). }\label{fig:dai3}
\end{figure}

Solutions to the Schr{\"o}der problem FEd($0$) are found,
in the usual method, by
composing a period-$Q$ function, $P$, with $a$. 
A typical example with $P$ chosen to be a cosine is shown in figure~\ref{fig:dai3}. The plotted streamfunction has many interesting nodal curves in addition to the nodal curve along the bottom topography $d(x)=1/|x|$ (black line). 
With the cosine $P$ there are 
elliptic nodal lines around the origin.

\section{Involutions, and a particularly simple family of solutions}\label{sec:Involutions}

Involutions are functions which when composed with themselves give the identity function:
$${\rm invol}({\rm invol}(x))= x $$
for all $x$ in the domain of the function.

It has already been noted, e.g.~\cite{ML95},  that
everywhere subcritical symmetric profiles lead to functional equations
$f(x)=f(T(x) ) $ where $T(x)=-{\rm invol}(x)$:
various examples are treated  in \S\ref{sec:subcritical}.
We do not know of any general method which is convenient to apply for all
equations of this type. 
If one simply changes the minus to a plus, we will see that the equation is extremely easy to solve.
\medskip

\begin{theorem}{\label{thm:invol1}
There are no solutions to the Abel functional equation, 
with $Q\ne{0}$
$$ a( {\rm invol}(x) ) - a(x) = Q $$}
\end{theorem}
\par\noindent{\it Proof.}
Suppose there were to be a solution to the Abel functional equation above, then
we also have
$$ a(x) - a( {\rm invol}(x) ) = a( {\rm invol}( {\rm invol}(x) ) ) - a( {\rm invol}(x) ) = Q$$
Adding the two preceding equations gives $0=2Q$ which contradicts the
assumption $Q\ne{0}$.
\medskip

Because of the preceding result, the approach --  using a solution of the
Abel equation to generate solutions to the
Schr{\"o}der equation by compositions with periodic functions --
fails here. However an alternative approach is available:
\smallskip

\begin{theorem}{ {\label{thm:invol2}}
Let $S$ be any symmetric function of two variables, meaning that
$S(u,v)=S(v,u)$ for all $u$, $v$.
Then the function $f(x)=S(x, {\rm invol}(x) ) $ solves the Schr\"{o}der equation
\begin{equation}
f({\rm invol}(x))=f(x) \quad{\rm with\ } {\rm invol} \quad{\rm an\ involution} .
\label{eq:involschroder}
\end{equation}
}
\end{theorem}
\par\noindent{\it Proof.}
$$f({\rm invol}(x))=S({\rm invol}(x),{\rm invol}({\rm invol}(x))) 
=S({\rm invol}(x),x)=S(x,{\rm invol}(x))=f(x) .$$

\medskip

For invol($x$) to correspond to a forward map $T$ we need to make sure that its domain is so that ${\rm invol}(x)>x$.

The entries in
the table below
indicate some flows associated with the involutions given.
We take $\nu=1$.
The entry $d$ is the solution of ${\rm invol}(x-d)=x+d$ .
There are many possibilities for $S$; our descriptions of the flow are for
$S(u,v)=u+v$.
(Any streamfunction $\sfun$ defined by the usual extension of $f$ is zero on $z=-d(x)$.)

\medskip


\noindent{}
\begin{tabular}{|| c | c |  l ||}
\hline
${\rm invol}(x)$& $d$ & 
Comments \\
\hline
$\frac{1}{x}$& $\sqrt{x^2-1}$ for $x<-1$& 
corner flow with a hyperbolic boundary\\
$\frac{x_0-x}{1+b x}$& $\sqrt{(x+\frac{1}{b})^2-\frac{1+ x_0 b}{b^2}}$&
further flows with hyperbolic $d$\\
$\sqrt{2b^2-x^2}$& $\sqrt{b^2-x^2}$&
$d$: portion of ellipse\\
${\rm PL}(x_0,m,x)$ with $m>1$& $\frac{(m+1)(x-x_0)}{m-1}$& 
piecewise linear $\sfun$ giving a corner\\
& & 
flow in a supercritical wedge\\
\hline
\end{tabular}

Some comments on the table above follow:\\
$\bullet$ Concerning the third entry in the table, we remark that
Barcilon's solution in a  circular quadrant with $\nu=1$ can be constructed using the discontinuous involution
${\rm sign}(x)\sqrt{1-x^2}$ and $f(x)=x^4 + {\rm invol}(x)^4$.\\
$\bullet$ In the fourth entry in the table, the piecewise linear involution PL is defined, with $m>1$, by
$$
{\rm PL}(x_0,m,x) = \frac{1}{2} \left(m-\frac{1}{m}\right) \left| x_0-x\right| +\frac{1}{2}
   \left(m+\frac{1}{m}\right) (x_0-x)+x_0
$$
There are several ways to generate the piecewise linear $\sfun$ corner flow.
One might take the symmetric function $S$ as
$S(u,v)=u+v$ or, alternatively, as $S(u,v)={\rm min}(u,v)$.
Let $\Gamma$ be the characteristic through $(x_0,0)$ extending downwards and to the right.
The flow has its streamlines parallel to $z=0$ in the triangle below the top boundary and above $\Gamma$
 and parallel to the bottom profile $z=-d(x)$ in the triangle above it and below $\Gamma$.
Taking $f=x+{\rm PL}(x_0,m,x)$  generates a similar corner flow. 
\smallskip

The corner flows, with no interior nodal lines, can be composed with other functions,
e.g. periodic functions, and then the $\sfun$ has nodal curves --
the flow exhibiting cells as in many of our earlier examples.

\medskip

Functions whose $k$-th iterate, $k\ge{2}$ is the identity are called
{\it involutions of order $k$}. 
The account above treats the case $k=2$, and it generalises.
For any $k\ge{2}$ there are no solutions to the 
involution Abel equations with $Q\ne{0}$.
Also, let $S$  be a function of $k$ arguments which is invariant
as one cycles through them, 
$$S(u_1,u_2,u_3,\ldots , u_k)= S(u_2,u_3\ldots, u_{k},u_1). $$
define
$$f(x)= S(x,{\rm invol}_k(x), {\rm invol}_k^{[2]}(x),\ldots{\rm invol}_k^{[k-1]}(x) ). $$
Then, for any $k\ge{2}$, $f$ 
solves FET($0$) when $T={\rm invol}_k$ is an involution of order $k$.
(Examples of $S$ include symmetric functions such as the sum of $k$ variables, etc..)

\section{Discussion}\label{sec:Discussion}

Solutions to the functional equations FEd($Q$) and FET($Q$) can be used to construct exact two-dimensional standing internal wave solutions. Several approaches for subcritical and (partly) supercritical domains making use of the functional equations are presented. There are others, e.g. the iterative methods due to Levy and others (see~\cite{Ku68}).
We believe that our exposition of the methods is satisfactory in the case of everywhere subcritical bottom profiles,
our \S{\ref{sec:Wunsch}} and \S\ref{sec:subcritical}: 
these are solutions where the `rays focus to the endpoints'.
For partly supercritical bottom profiles -- where the determination of the values of $\nu$ for which there are
solutions is also part of the problem -- our examples suggest that the functional equation approach may have value.
Our work on this in \S\ref{sec:partSuper} and \S\ref{sec:Involutions}
 is as much intended to publicise the problem as to present solutions.

The functional equations FEd($0$) and FET($0$) have been used in the past to construct exact internal wave solutions, and  \cite{Sa76} has also pointed out that one can associate solutions to FET($0$) with solutions to FET($Q\neq 0$). What is new is to link FET($Q\neq 0$) to Abel's functional equation and to make use for known properties and solutions of Abel's functional equation. Theorem~\ref{thm:A2} guarantees that for subcritical bottom topographies all solutions to FET($0$) are derived by applying the set of periodic function with period $1$ to any injective continuous solution of FET($1$). We are convinced that there is more to be elaborated, especially with the results on Abel's functional equation in \cite{Ku68}. 

We expect that functional equation techniques will prove useful for some other internal wave problems in which $z=0$ is a streamline.\\
1) One such situation concerns the generation of internal waves by
horizontal oscillations of a symmetric cylinder.
The usual formulation has the stream function $\psi_{\rm gen}$ nonzero on the cylinder:
 $\psi_{\rm gen}=-Uz$ on the cylinder $z=\pm{d(x)}$:
 see equation~(2.7) of~\cite{Hu97}.
 The \pde remains the wave equation as in our equations~(\ref{eq:pdeD}),
 but the boundary conditions, except for $\psi_{\rm gen}(x,0)=0$ are different.
 The representation of solutions as in equation~(\ref{eq:dalam}) with the boundary condition on
 the cylinder yields the functional equation
 $$ f_{\rm gen}(x-\frac{d(x)}{\nu})- f_{\rm gen}(x+\frac{d(x)}{\nu}) = -U d(x) .$$
 One solution of this is of the form $f_{\rm gen}(x)=c_{\rm gen} x$ with the constant $c_{\rm gen}=\nu U/2$.
 If $f$ solves the homogeneous equation FEd(0) then the general solution of the
 displayed equation immediately above is $f_{\rm gen}(x)=c_{\rm gen} x +f(x)$.
 The problem now requires complex-valued solutions of the functional equation
 with appropriate behaviour at infinity, a radiation boundary condition there.
 Several special cases have been investigated, and some solved by other techniques.
 \begin{itemize}
 \item
Elliptical cylinders with axes aligned with the coordinate axes are a particular case of the
more general treatment in ~\cite{Hu97}.
Here consider only the case
when $V=0$ in equation (3.42).
The $\sigma_\pm$ in~\cite{Hu97} is a multiple of our 
${x}\pm{z/\nu}$: see his equation (3.3).
Barcilon's (real) polynomial solutions correspond to blinking modes.
For the wave-generation problem of~\cite{Hu97} the complex-valued $f$ requires careful treatment of branch cuts in order that the
radiation conditions at infinity are satisfied.
\item An experimental treatment of a square cylinder is given in~\cite{Dal}.
\end{itemize}
2) Another instance where complex $f$, and radiation conditions,
are involved is the propagation, transmission and reflection of
monochromatic internal waves in a channel with a rigid upper lid 
$\{(x,0) |\, -\infty< x<\infty\}$ and an everywhere subcritical bottom, 
see~\cite{ML00,BH11}.
\smallskip

Further discussion is given in~\cite{Ke14}

\par\noindent
{\small{\bf Acknowledgements: } FB is grateful to Leo Maas and Gerard Sleijpen for comments on the manuscript, to the Institut for Marine and Atmospheric research Utrecht (IMAU) for providing a suitable workplace and to his fellow students at IMAU for creating an inspiring  working atmosphere.
GK is grateful to the  Centre for Water Research (CWR) at the University of Western Australia for the extended visit from
Curtin University during which this article was written.
An early version of this article is CWR reference 2672.}


\end{document}